\documentclass[11pt]{amsart}
\usepackage{geometry}                
\usepackage{graphicx}
\usepackage{amsmath, amsthm, amssymb, amscd, latexsym}
\usepackage{amssymb}

\numberwithin{equation}{section}
\newtheorem{theorem}{Theorem}[section]
\newtheorem{corollary}[theorem]{Corollary}
\newtheorem{lemma}[theorem]{Lemma}
\newtheorem{proposition}[theorem]{Proposition}
\newtheorem{conjecture}[theorem]{Conjecture}

\theoremstyle{definition}
\newtheorem{definition}[theorem]{Definition}
\newtheorem{remark}[theorem]{Remark}
\newtheorem{example}[theorem]{Example}

\usepackage{hyperref}
\usepackage{color}
\usepackage{cite}

\hypersetup{colorlinks=true,
     breaklinks=true,
     linkcolor=blue,    
     bookmarks=true,
      pageanchor=true}

\begin{document}

\title[On the Lehmer's totient problem on Number Fields]{On the Lehmer's totient problem on Number Fields}

\author[Konstantinos Smpokos]{Konstantinos Smpokos, Mathematician, PhD}

\email{k.smpokos@gmail.com }



\begin{abstract}
Lehmer's totient problem asks if there exists a composite number $d$ such that its totient divide $d-1$. In this article we generalize the Lehmer's totient problem in algebraic number fields. We introduce the notion of a Lehmer number. Lehmer numbers are defined to be the natural numbers which obey the Lehmer's problem in the ring of algebraic integers of a number field.
\end{abstract}

\maketitle

\section{Introduction}
Lehmer's totient problem is an open conjecture in Number Theory. It states that if for a natural number $d$,  $\phi(d) \mid d-1$ then $d$ is a prime number. Here $\phi$ is the Euler's totient function. In this article we will generalize Lehmer's totient problem on algebraic number fields. The original Lehmer's totient conjecture is equivalent to the generalized problem in $\mathbb{Q}$. Also, we define the concept of a Lehmer number over $K$. Then the generalized Lehmer's problem in $K$  is equivalent to say that every natural number is a Lehmer number over $K$. The generalized problem in a number field which is not $\mathbb{Q}$ is linked to the value of the Riemann zeta function in the degree of the extension of $K$ over $\mathbb{Q}$. The main reason we can not apply the same argument in $\mathbb{Q}$ is that the Riemann zeta function is infinite at 1. For the generalized problem in a number field $K$, we generalize the Euler's totient function $\phi$ and denote it by $\phi_K$. Of course, we have $\phi_{\mathbb{Q}} = \phi$. We consider the generalized problem only in number fields in which their ring of algebraic integers is a unique factorization domain. This is because in order to talk about the generalized Euler's totient function we need the concept of the greatest common divisor. In the generalized problem integers become algebraic integers of the number field.
\par
 In this article, we will define the concept of a realizable number over $K$. A natural number is defined to be realizable over $K$ if  every prime divisor of this number in $\mathbb{Z}$ is irreducible element in the ring of algebraic integers of $K$. We will show that for a number field which is not $\mathbb{Q}$, if a  natural number is realizable then it is a Lehmer number. Thus, the generalized problem is linked to primes in $\mathbb{Z}$ which are irreducible elements of the domain of algebraic integers of a number field. As a result, we discuss the Lehmer's problem in $\mathbb{Q}(i)$. In this article we define the notion of a normal number over $K$. We prove that $d$ is a Lehmer number over $\mathbb{Q}$ if and only if there exists a number field $K$ such $d$ is a normal Lehmer number over $K$. Also, we prove that  if $[K:\mathbb{Q}] \geq 2$, then $K$ is a realizable field if and only if $K$ is a Lehmer field and every prime in $\mathbb{Z}$ is a normal number over $K$. Also, we define the notion of a strongly Lehmer number. We have that in $\mathbb{Q}$ a natural number is strongly Lehmer number if and only if it  is a Lehmer number. However, in an arbitrary number field these two notions do not coincide in general. Also, we prove that if $[K:\mathbb{Q}] \geq 2$ then $K$ is a realizable field if and only if $K$ is a strongly Lehmer field. We prove that if $K,L$ are number fields such that the extensions $K/\mathbb{Q}$, $L/\mathbb{Q}$ are isomorphic, then $K$ and $L$ have the same realizable numbers. Finally, we will prove a result concerning the Lehmer's totient problem over $\mathbb{Q}$.

\section{Congruences in number fields}
To generalize the Euler's totient function to a number field we need to define congruences over the algebraic integers of the field. The definition of congruences is natural. We denote the domain of algebraic integers of the field $K$ by $O_K$. From now on we will assume that $O_K$ is a unique factorization domain. In a unique factorization domain the greatest common divisor of two algebraic integers is always defined. From now on if $x,y \in O_K$ we write $((x,y))$ for the greatest common divisor of $x$ and $y$.
\begin{definition}
Let $K$ be an algebraic number field. Suppose that we have $d \in \mathbb{N}$ and $x,y \in O_K$.
We define $x \equiv y \, (\mathsf{mod}\, d)$ iff $d \mid x-y$ in $O_K$.
\end{definition}
It is easy to see that $\equiv$ is an equivalence relation.
Let $d$ be a natural number. The set of equivalence classes of the relation $\equiv$ in $O_K$ is denoted by $\mathbb{Z}_d|_K$.
\begin{definition}
 If $x$,$y$ $\in O_K$ then in the set $\mathbb{Z}_d|_K$ we define $[x]_d + [y]_d = [x+y]_d$  and $ [x]_d [y]_d = [xy]_d$, where $[x]_d$ and $[y]_d$ are the equivalence classes of $x$ and $y$ respectively.
\end{definition}
Under these operations, it is easy to see that $\mathbb{Z}_d|_K$ is a ring.
We will prove later that $\mathbb{Z}_d|_K$ is a finite ring for every number field $K$.
Now we can define the unit group of integers modulo $d$ in $K$. The elements of this group are equivalence classes of algebraic integers of the number field.
\begin{definition}
Let $d$ be a natural number and $K$ be a number field.
The set $U(\mathbb{Z}_d)|_K$ is defined to be the set with elements all $[x]_d \in  \mathbb{Z}_d|_K$ such that $((x,d))=1$. The set $U(\mathbb{Z}_d)|_K$ is called the unit group of integers modulo $d$ in $K$.
\end{definition}
Like in $\mathbb{Q}$ we have that $U(\mathbb{Z}_d)|_K$ is a group under multiplication. It is evident to see that $U(\mathbb{Z}_d) = U(\mathbb{Z}_d)|_{\mathbb{Q}}$.
We have the following Proposition.
\begin{proposition} \label{prop1-}
Let $K$ be a number field and $d$ be a natural number. Then the ring $\mathbb{Z}_d|_K$ is a finite ring. Moreover, if $[K:\mathbb{Q}]=n$, then $|\mathbb{Z}_d|_K|= d^n$ and also $\mathbb{Z}_d|_K \cong (\mathbb{Z}_d)^n$ as additive groups, where $\mathbb{Z}_d$ is the usual ring of rational integers $\mathsf{mod}\, d$.
\begin{proof}
From \cite{Stewart} we have that the field $K$ has an integral basis. Then, if the integral basis is $\{w_1,..,w_n\}$ we have that  every algebraic integer $x$ can be written uniquely in the form $x= \sum_{i=1}^{n} k_i w_i$, where $k_i$ are rational integers. Now, since the numbers $k_i$ are rational integers we get $k_i \equiv \lambda_i \, (\mathsf{mod}\, d)$, for some $\lambda_i \in \{0,1..,d-1\}$.
Thus, $x= \sum_{i=1}^{n} k_i w_i \equiv \sum_{i=1}^{n} \lambda_i w_i , (\mathsf{mod}\, d)$. Now, the number of elements of the form $\sum_{i=1}^{n} \lambda_i w_i$ is finite, because $\lambda_i \in \{0,1,..,d-1\}$. Therefore, it is evident that the ring $\mathbb{Z}_d|_K$ is a finite ring.
For the second result we define the map $f: \mathbb{Z}_d|_K \rightarrow (\mathbb{Z}_d)^n$ with $f(\sum_{i=1}^{n}\lambda_i w_i) = (\lambda_1,..,\lambda_n)$.
Now we can see that this mapping is a group monomorphism under addition. Also, the number of the elements of $\mathbb{Z}_d|_K$ is $d^n$, since every number in $\mathbb{Z}_d|_K$ can be written uniquely in the form $\sum_{i=1}^{n}[\lambda_i]_d [w_i]_d$. Thus the rings  $\mathbb{Z}_d|_K$ and $(\mathbb{Z}_d)^n$ have the same number of elements. This gives that the mapping $f$ is also surjective. Thus, it is a group isomorphism under addition.
\end{proof}
\end{proposition}
Now, since we proved that integers modulo $d$ in $K$ form a finite ring, it is evident to see that $U(\mathbb{Z}_d)|_K$ is a finite group. This is because $U(\mathbb{Z}_d)|_K \subseteq \mathbb{Z}_d|_K$. Now we are ready to define the generalized totient function over a number field.
\begin{definition} \label{def1}
Let $K$ be a number field and $d \in \mathbb{N}$. The generalized totient function $\phi_K(d)$ is defined to be the order of the group $U(\mathbb{Z}_d)|_K$. That is, $\phi_K(d) = |U(\mathbb{Z}_d)|_K|$. It is evident that $\phi_{\mathbb{Q}} = \phi$.
\end{definition}
\begin{lemma} \label{lemma1}
Let $d \in \mathbb{N}$ and $K$ number field with $[K:\mathbb{Q}]=n$. Then $d$ is irreducible in $O_K$ if and only if $\phi_K(d) = d^n-1$.
\begin{proof}
Suppose $d$ is irreducible in $O_K$. Then $((w,d)) = 1$ for every $[w]_d \in \mathbb{Z}_d|_K$ with $[w]_d \neq [0]_d$. Thus, $U(\mathbb{Z}_d)|_K = (\mathbb{Z}_d|_K)^*$ which gives $\phi_K(d) = d^n-1$  by Definition \ref{def1} and  Proposition \ref{prop1-}.
Now, suppose $\phi_K(d) = d^n-1$. This gives that $U(\mathbb{Z}_d)|_K = (\mathbb{Z}_d|_K)^*$ since $U(\mathbb{Z}_d)|_K \subseteq  (\mathbb{Z}_d|_K)^*$.
This gives that $((w,d)) =1$ for every $w \in O_K$ with $[w]_d \neq [0]_d$. Therefore, it is evident that $d$ is irreducible in $O_K$.
\end{proof}
\end{lemma}
\begin{corollary}
$\mathbb{Z}_d|_K$ is a field if and only if $d$ is irreducible in $O_K$.
\end{corollary}
\section{The generalized totient function on Number Fields}
As we have seen, the generalized totient function $\phi_K$ is the order of the group $U(\mathbb{Z}_d)|_K$.
In this section we will prove some properties of $\phi_K$ similar to the properties of the Euler's totient function $\phi$.
\begin{proposition} \label{prop2-}
Let $K$ be a number field and $m.n \in \mathbb{N}$ relatively prime. Then we have $\phi_K(mn) = \phi_K(m) \phi_K(n)$. That is, $\phi_K$ is a multiplicative arithmetic function for every number field $K$. Moreover, we have $\mathbb{Z}_{mn}|_K \cong \mathbb{Z}_m|_K \times \mathbb{Z}_n|_K$ as rings and $U(\mathbb{Z}_{mn})|_K \cong U(\mathbb{Z}_m)|_K \times U(\mathbb{Z}_n)|_K$ as groups.
\begin{proof}
The proof idea comes from \cite{Stewart}.
Let a mapping $g: O_K \rightarrow \mathbb{Z}_m|_K \times \mathbb{Z}_n|_K$, with $g(x)= ([x_1]_m,[x_2]_n)$, where $x \equiv x_1 \, (\mathsf{mod}\, m)$ and $x \equiv x_2 \, (\mathsf{mod}\, n)$. We can see that this map is a well defined ring homomorphism. Let $x \in Kerg$. Then , $m \mid x_1$ and $n \mid x_2$. Thus, $m \mid x$ and $n \mid x$, by the definition of $x_1$ and $x_2$. Thus since $m,n$ are relatively prime we have $mn \mid x$. That is $x \in mn O_K$. This shows that $Kerg = mn O_K$. Therefore it is evident that $O_K/ mn O_K \cong img$. Thus, $\mathbb{Z}_{mn}|_K \cong img$ and also we have by Proposition \ref{prop1-} that $|\mathbb{Z}_{mn}|_K| = (mn)^{[K:\mathbb{Q}]} = m^{[K:\mathbb{Q}]} n^{[K:\mathbb{Q}]} = |\mathbb{Z}_m|_K \times \mathbb{Z}_n|_K|$. Thus, $img = \mathbb{Z}_m|_K \times \mathbb{Z}_n|_K$ and this shows that  $\mathbb{Z}_{mn}|_K \cong \mathbb{Z}_m|_K \times \mathbb{Z}_n|_K$ as rings.
Now, we can see that $U(\mathbb{Z}_{mn})|_K$ is the group of units of $\mathbb{Z}_{mn}|_K$ and  $U(\mathbb{Z}_m)|_K \times U(\mathbb{Z}_n)|_K$ is the group of units of $\mathbb{Z}_m|_K \times \mathbb{Z}_n|_K$. Therefore, it is evident that $U(\mathbb{Z}_{mn})|_K \cong U(\mathbb{Z}_m)|_K \times U(\mathbb{Z}_n)|_K$ as groups and this also shows that $\phi_K(mn) = \phi_K(m) \phi_K(n)$ and the proof is complete.
\end{proof}
\end{proposition}
The next Proposition relates the Euler's totient function with the generalized totient function over a number field $K$.
\begin{proposition}
Let a number field $K$ and $d \in \mathbb{N}$. The unit group of rational integers modulo $d$ can be embedded to $U(\mathbb{Z}_d)|_K$.
Thus, we have $\phi(d) \mid \phi_K(d)$ for every $d$ and every $K$ number field.
\begin{proof}
Let $j: U(\mathbb{Z}_d) \rightarrow U(\mathbb{Z}_d)|_K$ with $j([x]_d) = [x]_{d,K}$.
Now, if $x \in \mathbb{Z}$ and $x$ and $d$ are relatively prime in $\mathbb{Z}$, then they are relatively prime in $O_K$, since $O_K$ is a principal ideal domain. This shows that $j$ is a well defined group homomorphism.
Assume $[x]_{d,K} = [1]_{d,K}$ for some $x \in \mathbb{Z}$. Then $d \mid x-1$ in $O_K$.Thus, there exists $y \in O_K$ such that $x-1 = dy$. Thus, we have that $y$ is a rational integer and this gives that $d \mid x-1$ in $\mathbb{Z}$. Therefore $j$ is a group monomorphism and the proof is complete.
\end{proof}
\end{proposition}
\section{The Lehmer's totient problem and its generalization}
Now we will state the Lehmer's totient problem and we will generalize this open problem in algebraic number fields.
\begin{conjecture} (Lehmer's totient problem)
\\
Assume $d \in \mathbb{N}$. Then if $\phi(d) \mid d-1$ then $d$ is a prime number.
\end{conjecture}
We have already showed that a natural number $d$ is irreducible in $O_K$ if and only if $\phi_K(d) = d^n-1$, where $n= [K:\mathbb{Q}]$. Also, we have that $d^n-1$ is the maximal value that $\phi_K(d)$ can take. This is similar to $\mathbb{Q}$, where the maximal value $\phi(d)$ can take is $d-1$, which is attained if and only if $d$ is a prime number.
 It is easy to show that if for a natural number $d$ is true that $\phi(d) \mid d-1$, then we have that $d$ is squarefree. Therefore, from now on we are concerned only with squarefree numbers.
This leads us to the following definition of a Lehmer number.
\begin{definition} \label{def2}
Let $d \in \mathbb{N}$  and $K$ be a number field with $[K:\mathbb{Q}]=n$. We define that $d$ is a Lehmer number over $K$ if and only if the following statements are equivalent.
\\
1. $\phi_K(d) \mid d^n-1$.
\\
2. d is an irreducible element of $O_K$.
\end{definition}
\begin{remark} \label{remark1}
In the light of Definition \ref{def2} we can easily see that the Lehmer's totient conjecture is true if and only if every natural number is a Lehmer number over $\mathbb{Q}$.
This article examines the case where $[K:\mathbb{Q}] \geq 2$. There may be a connection from the number fields $K$ with $[K:\mathbb{Q}] \geq 2$ with the case $K=\mathbb{Q}$, but we are not aware of it in this moment.
\end{remark}
\begin{definition} (Generalized Lehmer's totient problem)
\\
Let $K$ be a number field. The field $K$ satisfies the generalized Lehmer's totient problem if and only if every natural number is a Lehmer number over $K$.
Such fields are called Lehmer fields.
\end{definition}
\begin{definition} \label{def3}
Let $d$ be a natural number. We say that $d$ is realizable over $K$ if and only if all prime divisors of $d$ in $\mathbb{Z}$ are irreducible elements of $O_K$.
\begin{remark} \label{remark2}
It is evident that all natural numbers are realizable over $\mathbb{Q}$. We will show that if $[K:\mathbb{Q}] \geq 2$ then every realizable number is a Lehmer number. To avoid confusion we should say that realizable numbers over $\mathbb{Q}$ are not necessairily Lehmer numbers over $\mathbb{Q}$. Note the $[K:\mathbb{Q}] \geq 2$ condition.
\end{remark}
\begin{definition}
Let $K$ be a number field. We say that $K$ is realizable iff all  natural numbers are realizable over $K$.
\end{definition}
\end{definition}
Now we find a characterization of realizable numbers relevant to the generalized totient function over a number field.
\begin{proposition} \label{prop12}
Let $d \in \mathbb{N}$, squarefree and $K$ be a number field with $[K:\mathbb{Q}] = n$. Then $d$ is realizable over $K$ if and only if $\phi_K(d) = \prod_{p \mid d} (p^n-1)$, where the product is over all prime divisors of $d$ in $\mathbb{Z}$.
\begin{proof}
Assume that $d$ is realizable over $K$. By Proposition \ref{prop2-} we have $\phi_K(d) = \prod_{p \mid d} \phi_K(p)$. Then since $d$ is realizable we have that each prime divisor of $d$ in $\mathbb{Z}$ is an irreducible element of $O_K$. Thus, if $p \mid d$ we have $\phi_K(p) = p^n-1$, by Lemma \ref{lemma1}. Therefore, $\phi_K(d) = \prod_{p \mid d}(p^n-1)$.
Now, since for a prime number $p$ in $\mathbb{Z}$ we have $\phi_K(p) \leq p^n-1$ with equality if and only if $p$ is irreducible in $O_K$, we get that if $d$ is not realizable then $\phi_K(d) < \prod_{p \mid d} (p^n-1)$ and the proof is complete.
\end{proof}
\end{proposition}
For the following theorem we need the well known Riemann zeta function defined on real numbers $s>1$ discussed in \cite{Apost}. We have that if $s>1$, $\zeta(s) := \sum_{n=1}^{\infty} \frac{1}{n^s}$. Riemann's zeta function can be written as an absolutely convergent Euler product as $\zeta(s) = \prod_{p: prime}\frac{1}{1- p^{-s}}$.
\begin{theorem} \label{theorem1}
Let $d \in \mathbb{N}$ and $K$ be a number field with $[K:\mathbb{Q}] \geq 2$.
We have that if $d$ is realizable over $K$ then $d$ is a Lehmer number over $K$.
\begin{proof}
From the assumptions and Proposition \ref{prop12} we have $\phi_K(d) = \prod_{p \mid d}(p^n - 1)$. 
Therefore, $\frac{d^n - 1}{\phi_K(d)} =$$ \frac{(\prod_{p \mid d} p)^n -1}{\prod_{p \mid d}(p^n - 1)} \leq \prod_{p \mid d} \frac{1}{1- p^{-n}}$ $\leq \prod_{p: prime} \frac{1}{1- p^{-n}} = \zeta(n)$. Thus, $\frac{d^n-1}{\phi_K(d)} \leq \zeta(n) $.
Assume that $\phi_K(d) \mid d^n-1$. Assume that $d$ is not irreducible element in $O_K$. Thus, $\frac{d^n-1}{\phi_K(d)} \geq 2$.
This gives that $2 \leq \zeta(n)$, which is a contradiction since $\zeta$ is stricly decreasing on $\{w \in \mathbb{R}: w>1\}$ and $n \geq 2$.
Now the proof is complete.
\end{proof}
\end{theorem}
\begin{corollary}
If a number field $K \neq \mathbb{Q}$ is realizable then it is a Lehmer field.
\end{corollary}
\begin{example}
We take $K= \mathbb{Q}(i)$. Then we have by \cite{Rings} that a prime in $\mathbb{Z}$ is irreducible in $O_K$ if and only if $p \equiv 3 \, (\mathsf{mod}\, 4)$.
This gives that in $\mathbb{Q}(i)$ all numbers in which every prime divisor is $\equiv 3 \, (\mathsf{mod}\, 4)$ are Lehmer numbers.
\end{example}
\begin{definition} \label{def6}
Let $K$ be a number field with $[K:\mathbb{Q}] = n$ and $d$ be a natural number. We call the number $d$ normal in $K$ iff $\frac{\phi_K(d)}{\phi(d)} \mid \frac{d^n-1}{d-1}$.
We say that $K$ is a normal field iff every natural number is normal in $K$.
\end{definition}
Using this definition we prove the following Theorem.
\begin{theorem} \label{thm2}
Let $d$ be a natural number. Then, we have that $d$ is a Lehmer number over $\mathbb{Q}$ if and only if there exists a number field $K$ such that $d$ is Lehmer and normal over $K$.
\begin{proof}
The only if part is obvious, since every natural number is normal over $\mathbb{Q}$.
Now, assume that there exists a number field $K$ such that $d$ is Lehmer and normal over $K$.
Assume that $\phi(d) \mid d-1$. Then by the normality of $d$ we have that $(d-1)\phi_K(d) \mid (d^n-1) \phi(d)$. By assumption we have $\phi_K(d) \phi(d) \mid (d-1) \phi_K(d) \mid (d^n -1) \phi(d)$ which gives that $\phi_K(d) \mid d^n-1$. Therefore, we get that $d$ is irreducible in $O_K$, since $d$ is a Lehmer number over $K$. Thus, $d$ must be prime and the proof is complete.
\end{proof}
\end{theorem}
In the next Lemma we characterize the primes which are normal elements over a number field $K$.
\begin{lemma} \label{lem1}
Assume $[K:\mathbb{Q}] =  n$.
\\
If $p$ is prime in $\mathbb{Z}$, then $p$ is a normal number over $K$ if and only if $\phi_K(p) \mid p^n-1$.
\begin{proof}
The proof is straightforward.
By definition \ref{def6}, we have $p$ is normal over $K$ if and only if $\frac{\phi_K(p)}{\phi(p)} \mid \frac{p^n-1}{p-1}$ which is equivalent to $\phi_K(p) \mid p^n-1$, since $p$ is a prime number.
\end{proof}
\end{lemma}
Now, we will characterize the realizable Number fields.
\begin{proposition} \label{prop1}
Let $K$ be a number field with $[K:\mathbb{Q}] =n \geq 2$. Then, the following are equivalent
\\
1. K is a realizable field
\\
2. K is a Lehmer field and every prime in $\mathbb{Z}$ is normal in $K$.
\begin{proof}
For the "only if" part. Assume that $K$ is a realizable field. By Theorem \ref{theorem1} we have that $K$ is a Lehmer field. Also, if $p$ is a prime number in $\mathbb{Z}$, then $p$ is an irreducible element in $O_K$ by the realizability of $K$. Therefore, $\phi_K(p) = p^n -1$ which gives that $p$ is a normal number in $K$.
Thus, the field $K$ is a Lehmer field and every prime number in $\mathbb{Z}$ is normal in $K$.
\\
For the "if" part, we have that if $p$ is a prime number in $\mathbb{Z}$, then $p$ is a normal element in $K$. Therefore, by Lemma \ref{lem1} we get $\phi_K(p) \mid p^n-1$. Also $p$ is a Lehmer number over $K$. Therefore, we get that $p$ is an irreducible element of $O_K$. Therefore, $p$ is realizable over $K$. Thus, $K$ is a realizable field and the proof is complete.
\end{proof}
\end{proposition}
Now we will define the notion of a strongly Lehmer number over a number field.
\begin{definition} \label{def7}
Let $d$ be natural number and $K$ be a number field with $[K:\mathbb{Q}]=n$. We say that $d$ is a strongly Lehmer number over $K$ if the following are equivalent.
\\
1. $\phi_K(d) \mid d^n-1$.
\\
2. $d$ is an irreducible number in $O_K$.
\\
3. $d$ is a prime number in $\mathbb{Z}$.
\end{definition}
It is evident that strongly Lehmer numbers over $K$ are Lehmer numbers over $K$ for every number field $K$.
\begin{theorem} \label{thm3}
Let $K$ be a number field with $[K:\mathbb{Q}] \geq 2$. Then $K$ is a realizable field if and only if  $K$ is a strongly Lehmer field.
\begin{proof}
For the "only if" part we have the following. Let $d$ be natural number. By Theorem \ref{theorem1} we get that $d$ is a Lehmer number over $K$. Also, $d$ is irreducible number over $O_K$ if and only if $d$ is a prime number in $\mathbb{Z}$ by the realizability of $K$. Thus, by Definition \ref{def7} we get that $d$ is a strongly Lehmer number over $K$.
\\
For the "if" part, let $p$ be a prime number in $\mathbb{Z}$. Then $p$ is a strongly Lehmer number over $K$. Thus, $p$ is a prime number if and only if $p$ is an irreducible element of $O_K$.
This proves that $p$ is realizable over $K$. Thus, $K$ is a realizable field.
\end{proof}
\end{theorem}
By \cite{galois} we have the following definition.
\begin{definition} \label{def8}
Let two field extensions $K/L$ and $M/N$ of number fields. Then, an isomorphism between these extensions is a pair $(\lambda,\mu)$ of field isomorphisms $\lambda:L \rightarrow N$, $\mu: K \rightarrow M$ such that $\lambda(x) = \mu(x)$ for every $x \in L$.
\end{definition}
\begin{corollary} \label{cor2}
Let $K,L$ be number fields. Then the extensions $K/\mathbb{Q}$ and $L/\mathbb{Q}$ are isomorphic iff there exists $\mu : K \rightarrow L$ field isomorphism such that $\mu(q)=q$ for every $q \in \mathbb{Q}$.
\end{corollary}
\begin{theorem} \label{thm4}
If the extensions $K/\mathbb{Q}$, $L/\mathbb{Q}$ are isomorphic and $d$ be a natural number, then $d$ is realizable over $K$ if and only if $d$ is realizable over $L$.
\begin{proof}
By Corollary \ref{cor2} we get that there exists $\mu : K \rightarrow L$ field isomorphism such that $\mu(q)=q$ for every $q \in \mathbb{Q}$.
Let $d$ be a natural number. Let $p$ prime with $p \mid d$. Assume that $d$ is a realizable number over $L$. Assume $p = \zeta_1 \zeta_2$, where $\zeta_1,\zeta_2 \in O_K$. Then we get $\mu(p)=p=\mu(\zeta_1) \mu(\zeta_2)$ in $O_L$. Thus, $\mu(\zeta_1)$ is a unit or $\mu(\zeta_2)$ is a unit, which gives that $\zeta_1$ is a unit or $\zeta_2$ is a unit. Thus, $p$ is realizable over $K$. Therefore $d$ is realizable over $K$.
By taking $\mu^{-1}: L \rightarrow K$ and by the same proof, we get that $p$ is realizable over $K$ only if $p$ is realizable over $L$. This completes the proof.
\end{proof}
\end{theorem}
\begin{corollary} \label{cor3}
Let $K,L$ be number fields such that $K/\mathbb{Q} \cong L/\mathbb{Q}$.
Then $K$ is a strongly Lehmer field if and only if $L$ is a strongly Lehmer field.
\end{corollary}
Now we will state a result from \cite{galois}, which we will use in the next Theorem.
\begin{proposition} \label{prop2}
Suppose that $K$ and $L$ are subfields of $\mathbb{C}$ and $i: K \rightarrow L$ is an isomorphism. Let $K(a),L(b)$ be simple algebraic extensions of $K$ and $L$ respectively, such that $a$ has minimal polynomial $m_a(t)$ over $K$ and $b$ has minimal polynomial $m_b(t)$ over $L$. Suppose $m_b(t) = i(m_a(t))$. Then there exists an isomorphism $j:K(a) \rightarrow L(b)$ such that $j|K = i$ and $j(a)=b$.
\end{proposition}
\begin{theorem} \label{thm5}
Let $a,b \in \mathbb{C}$. Assume that there exists an irreducible polynomial $p(t)$ in $\mathbb{Q}[t]$ such that $p(a)=p(b)=0$. Then a natural number is realizable over $\mathbb{Q}(a)$ if and only if it  is realizable over $\mathbb{Q}(b)$.
\begin{proof}
If we take $K=L=\mathbb{Q}$ and $i$ to be the identity map in Proposition \ref{prop2},  then by Thorem \ref{thm4} and Corollary \ref{cor2} we get the result.
\end{proof}
\end{theorem}
\section{Other Results}
In this section we will prove a result concerning the Lehmer's totient problem over $\mathbb{Q}$.
Here we will deal with squarefree numbers. The reason for this is that a counterexample to the Lehmer's totient problem is necessarily squarefree.
\begin{theorem} \label{thm6}
Let $w$ be a squarefree number and $l \in \mathbb{Q}$ such that $l < \frac{w}{\phi(w)}$. Then, there exist finitely many $d \in \mathbb{N}$ with $d$ squarefree such that $w \mid d$ and $\frac{d-1}{\phi(d)} = l$.
\begin{proof}
We have $\frac{d-1}{\phi(d)}$$ = \frac{(\prod_{p \mid d}p) -1}{\prod_{p \mid d}(p-1)}$ $= \prod_{p \mid d}(\frac{1}{1 - p^{-1}}) - \frac{1}{\phi(d)}$.
\\
For the sake of contradiction, we assume that there exist infinitely many $d$ squarefree with $w \mid d$ and $\frac{d-1}{\phi(d)} = l$. Then we can find a sequence of natural numbers with $\lim_{n \rightarrow \infty}d_n = \infty$ such that $d_n$ has this propetry.
Thus, since $\phi(d_n) \rightarrow \infty$ as $n \rightarrow \infty$, we get $\prod_{p \mid d_n}(\frac{1}{1 - p^{-1}}) \rightarrow l$, as $n \rightarrow \infty$.
Now since $w \mid d_n$, we get $\prod_{p \mid d_n}(\frac{1}{1 - p^{-1}}) \geq \prod_{p \mid w}(\frac{1}{1 - p^{-1}})$$= \frac{w}{\phi(w)}$.
By taking limits, we get $l \geq \frac{w}{\phi(w)}$, which is a contradiction by assumption.
\end{proof}
\end{theorem}
\

\end{document}